\documentclass[10pt,reqno]{amsart}

\usepackage[latin1]{inputenc}
\usepackage{amsmath}
\usepackage{amssymb}
\usepackage{hyperref}
\usepackage{stackrel}
\usepackage{xcolor}

\newtheorem{theorem}{Theorem}[section]
\newtheorem{lemma}[theorem]{Lemma}
\newtheorem{example}[theorem]{Example}
\newtheorem{definition}[theorem]{Definition}

\newtheorem{proposition}[theorem]{Proposition}
\newtheorem{remark}[theorem]{Remark}
\newtheorem{corollary}[theorem]{Corollary}

\begin{document}

\title{Sequence of Induced Hausdorff metrics on Lie groups}

\author{Norbil Cordova}
\address{Department of Mathematics, State University of Maring\'a,
87020-900, Maring\'a, PR, Brazil \\ email: nlcneyra2@uem.br}

\author{Ryuichi Fukuoka}
\address{Department of Mathematics, State University of Maring\'a,
87020-900, Maring\'a, PR, Brazil \\ email: rfukuoka@uem.br}

\author{Eduardo de A. Neves}
\address{Department of Mathematics, State University of Maring\'a,
87020-900, Maring\'a, PR, Brazil \\ email: eaneves@uem.br}

\date{November 9th, 2017}

\begin{abstract}
Let $\varphi: G \times (M,d) \rightarrow (M,d)$ be a left action of a Lie group on a differentiable manifold endowed with a metric $d$, which is compatible with its topology.
Let $X$ be a compact subset of $M$. Then the isotropy subgroup of $X$ is defined as $H_X:=\{g\in G; gX=X\}$ and it is closed in $G$.
The induced Hausdorff metric is a metric on the left coset manifold $G/H_X$ defined as $d_X(gH_X,hH_X)=d_H(gX,hX)$, where $d_H$ is the Hausdorff distance in $M$.
Suppose that $\varphi$ is transitive and that there exist $p\in M$ such that $H_X=H_p$.
Then $gH_X \mapsto gp$ is a diffeomorphism that identifies $G/H_X$ and $M$.
In this work we define a discrete dynamical system of metrics on $M$.
Let $d^1=\hat d_X$, where $\hat d_X$ stands for the intrinsic metric induced by $d_X$. 
We can iterate the process on $\varphi: G \times (M\equiv G/H_X,d^1)\rightarrow (M\equiv G/H_X,d^1)$, in order to get $d^2, d^3$ and so on.
We study the particular case where $M=G$,  $\varphi: G\times (G,d) \rightarrow (G,d)$ is the usual product, $d$ is bounded above by a right invariant intrinsic metric on $G$ and $X$ is a finite subset of $G$ containing the identity element.
We prove that $d^i$ converges pointwise to a metric $d^\infty$. 
In addition, if $d$ is complete and the semigroup generated by $X$ is dense in $G$, then $d^\infty$ is the distance function of a right invariant $C^0$-Carnot-Carath\'eodory-Finsler metric. 
The case where $d^\infty$ is $C^0$-Finsler is studied in detail.
\end{abstract}



\maketitle

\section{Introduction}
\label{introducao}

Before talking about the subject of this paper, we present two topics related to this work.
We don't use the second topic here. It will only illustrate our main result (Theorem \ref{teorema principal}).

The first topic is invariant metrics on homogeneous spaces. Let $M$ be a connected differentiable manifold endowed with a completely nonholonomic distribution $\mathcal D$. 
An absolutely continuous path $\gamma$ in $M$ is horizontal if $\gamma^\prime(t) \in \mathcal D$ for almost every $t$.
Chow-Rashevskii theorem states that every pair of points in $M$ can be connected by a  horizontal curve (see \cite{Chow}, \cite{Montgomery}, \cite{Rashevskii}).
We endow each subspace $\mathcal D_x$ of $\mathcal D$ with a norm $F(x,\cdot)$ such that $x \mapsto F(x,\cdot)$ is continuous, that is, given a horizontal continuous vector field $Y$ on $M$, then $x \mapsto F(x,Y(x))$ is a continuous map (see \cite{Berestovskii1}).
Now we proceed as in the definition of sub-Riemannian geometry: 
The $C^0$-Carnot-Carath\'eodory-Finsler metric on $M$ is given by
\begin{equation}
\label{equacao carnot caratheodory}
d_c(x,y)=\inf_{\gamma \in \mathcal H_{x,y}}\int_I F(\gamma(t),\gamma^\prime(t))dt,
\end{equation}
where $\mathcal H_{x,y}$ is the set of horizontal curves connecting $x$ and $y$ (see \cite{Berestovskii1}, \cite{Berestovskii2}).

The following theorem due to Berestovskii is important for this work and states that if an intrinsic metric is homogeneous, then it has the tendency to gain extra regularity.

\begin{theorem}
\label{Berestovskii theorem}[Berestovskii, \cite{Berestovskii2}, Theorem 3]
If $(M,d)$ is locally compact, locally contractible homogeneous space, endowed with an invariant intrinsic metric $d$, then $(M,d)$ is isometric to the quotient space $G/H$ of a Lie group $G$ over a compact subgroup $H$ of it, endowed with a $G$-invariant $C^0$-Carnot-Carath\'eodory-Finsler metric.
\end{theorem}

The second topic is geometric flows in Riemannian geometry. 
There are several of them such as the mean curvature flow (for hypersurfaces in Riemannian manifolds), the Ricci-Hamilton flow, etc. 
These flows, in some specific situations, converges to a more ``homogeneous'' Riemannian metric (eventually after some normalization of the volume).
For instance, the normalized Ricci-Hamilton flow converges to a Riemannian metric of constant sectional curvature for $3$-dimensional closed Riemannian manifolds with strictly positive Ricci curvature (see \cite{Hamilton1}).
It is an example of a dynamical system of metrics converging to a ``well behaved'' metric.

Now we give an outline of this work.

Let $\varphi: G \times M \rightarrow M$ be a left action of a Lie group on a differentiable manifold endowed with a metric $d$ (distance function) compatible with the topology of $M$.
As usual, we denote $gx:=\varphi(g,x)$. 
Let $X$ be a compact subset of $M$.
The isotropy subgroup of $X$ is a closed subgroup $H_X$ of $G$ given by
\[
H_X=\{g\in G; gX=X\}.
\]
Then there exist a unique differentiable structure on $G/H_X$, compatible with the quotient topology, such that the natural action $\phi: G \times G/H_X \rightarrow G/H_X$, given by $\phi(g,hH_X)=ghH_X$, is smooth.
In \cite{Benetti} and \cite{BenettiFukuoka}, the authors define a metric $d_X$ on $G/H_X$, which is called induced Haudorff metric, by 
\[
d_X(gH_X,hH_X)=d_H(gX,hX),
\]
where $d_H$ is the Hausdorff distance in $(M,d)$ and $gH_X$ denote the left coset of $g$ in $G/H_X$.
Suppose that $\varphi$ is transitive and that there exist $p\in M$ such that $H_p=H_X$.
Then the map $\eta: G/H_X \rightarrow M$ given by $gH_X \mapsto g.p$ is a diffeomorphism such that $\varphi(g,\eta(hH_X))= \eta(\phi(g,hH_X))$, that is, $\varphi$ can be identified with $\phi$.
We use the induced Hausdorff metric in order to create a discrete dynamical system of metrics on $M$.
Define $d^1=\hat d_X$, where the hat denotes the intrinsic metric induced by $d_X$.
We can iterate this process in $\varphi:G\times (M,d^1)\rightarrow (M,d^1)$ in order to obtain $d^2$, where $(G/H_X,d^1)$ is identified with $(M,d^1)$ via $\eta$.
Through this iteration, we can define a sequence of metrics $d,d^1,d^2, \ldots$ on $M$, which we call {\it sequence of induced Hausdorff metrics}.

In this work we study the sequence of induced Hausdorff metrics in the following case:
$M$ is the Lie group $G$ itself endowed with a metric $d$ that is bounded above by a right invariant intrinsic metric $\mathbf d$, $\varphi: G \times (G,d) \rightarrow (G,d)$ is the product of $G$ and
$X$ is a finite subset of $G$ containing the identity element $e$ (It is implicit here that $H_X=\{e\}$).
We prove that $d^i$ converges pointwise to a metric $d^\infty$.
Moreover if the semigroup $S_X$ generated by $X$ is dense in $G$ and $d$ is complete, then $d^\infty$ is the distance function of a right invariant $C^0$-Carnot-Carath\'eodory-Finsler metric.
We give a necessary and sufficient condition in order to $d^\infty$ be $C^0$-Finsler. 
In this case, an explicit formula for the Finsler metric is obtained.

Of course the comparison of our work with the Riemannian geometric flows can't be taken so literally. 
It is given here only to illustrate the dynamics of the sequence of induced Hausdorff metrics, which in some cases converges to a more ``well behaved'' metric. 
One important restriction is that it is defined only on left coset manifolds. On the other hand the metric $d$ is much more general than Riemannian metrics and the induced Hausdorff distance is the basic tool in order to overcome the lack of differentiability.

This work is organized as follows. 
In Section \ref{preliminares} we fix notations and we present definitions and results that are necessary for this work. 
In Section \ref{sequencia de metricas} we study some properties of the sequence of induced Hausdorff metrics.
In particular, we prove that it is increasing, and if $d$ is bounded above by a right invariant intrinsic metric, then it converges pointwise to a metric.
This limit metric is intrinsic if $d$ is complete.
In Section \ref{existencia sx denso}, we prove that every connected Lie group $G$ admits a finite subset $X \ni e$ such that $H_X=\{e\}$ and $\bar S_X = G$.
In Section \ref{Finlser fomula}, we study the case $\bar S_X=G$, where $d$ is complete and bounded above by a intrinsic right invariant metric. 
We prove that the sequence of induced Hausdorff metrics converges to a right invariant $C^0$-Carnot-Carath\'eodory-Finsler metric. 
The case where $d^\infty$ is $C^0$-Finsler is also studied here, as explained before.
In Section \ref{secao exemplos}, we give some additional examples in order to illustrate better this work.
Finally Section \ref{secao final} is devoted to final remarks.

The authors would like to thank Professor Luiz A. B. San Martin for some valuable suggestions.

\section{Preliminaries}
\label{preliminares}

In this section we present some definitions and results that are used in this work. They can be found in \cite{Burago}, \cite{BenettiFukuoka}, \cite{Helgason}, \cite{Jurdjevic},  \cite{KobayashiNomizu1}, \cite{KobayashiNomizu2} and \cite{Warner}.
For the sake of clearness, we usually don't give the definitions and results in their most general case.

Let $G$ be a group and $(M,d)$ be a metric space. 
Consider a left action $\varphi:G\times M \rightarrow M$ of $G$ on $M$.
Then $ex=x$ for every $x\in M$ and $(gh)x=g(hx)$ for every $(g,h,x)\in G\times G\times M$. 
Every $\varphi_g:=\varphi(g,\cdot)$ is a bijection. 
We say that a left action $\varphi:G \times M \rightarrow M$ is an action by isometries if every $\varphi_g$ is an isometry. 
Analogously we say that $\varphi$ is an action by homeomorphism if $\varphi_g$ is a homeomorphism for every $g\in G$. 

Let $(M,d)$ be a metric space. 
We denote the open ball with center $p$ and radius $r$ in $(M,d)$ by $B_d(p,r)$. 
The closed ball is denoted by $B_d[p,r]$. 
The topology induced by $d$ is denoted by $\tau_d$.
The closure of a subset $A$ in $(M,d)$ is denoted by $\bar A$. When more than one metric or topology are involved, for instance we have metrics $d$ and $\rho$ and a topology $\tau$ on $M$, we use terms like $d$-neighborhood, $\tau$-open subset, $\rho$-compact, etc.

Let $A,B$ be compact subsets of $(M,d)$. The Hausdorff distance between $A,B$ is given by
\[
d_H(A,B)=\max\left\{\sup_{x\in A} \inf_{y \in B} d(x,y), \sup_{y\in B} \inf_{x \in A} d(x,y) \right\}.
\]
It is well known that $d_H$ is a metric on the family of compact subsets of $M$.

Let $\varphi:G \times M \rightarrow M$ be a left action by homeomorphisms of a group $G$ on a metric space $(M,d)$.
If $X\subset M$ is a subset, then the isotropy subgroup of $G$ with respect to $X$ is defined by $H_X=\{g\in G;gX=X\}$.
Suppose that $X$ is a compact subset in $M$.
In Proposition 2.1 of \cite{BenettiFukuoka} (and in \cite{Benetti}), we define the induced Hausdorff metric on the left coset space $G/H_X$ as $d_X(gH_X,hH_X)=d_H(gX,hX)$, where $d_H$ is the Hausdorff distance in $M$.

A partition $\mathcal{P}$ of an interval $[a,b]$ is a subset $\{ t_0,\ldots,t_{n_{\mathcal P}}\} \subset [a,b]$ such that $a=t_0 < t_1 < \ldots < t_{n_{\mathcal P}} = b$. 
The norm of $\mathcal P$ is defined as $\vert \mathcal P\vert=\max_{i=1,\ldots,n_{\mathcal P}} \vert t_i - t_{i-1}\vert$.
The length of a path $\gamma:[a,b] \rightarrow M$ on a metric space $(M,d)$ is given by
\[
\ell_d(\gamma)=\sup_{\mathcal{P}}\sum_{i=1}^{n_{\mathcal{P}}}d(\gamma(t_i),\gamma(t_{i-1})).
\]
It is well known that for every $\varepsilon >0$, there exist a $\delta >0$ such that
\[
\ell_d(\gamma)\leq \sum_{i=1}^{n_{\mathcal P}} d(\gamma(t_i),\gamma(t_{i-1}))+\varepsilon
\]
for every partition $\mathcal P$ such that $\vert \mathcal{P}\vert < \delta$ (see \cite{Burago}).

Given a metric space $(M,d)$ we can define the extended metric (the distance can be $\infty$)
\[
\hat d(x,y)=\inf_{\gamma \in \mathcal{C}_{x,y}}\ell_d(\gamma),
\]
on $M$, where $\mathcal{C}_{x,y}$ is the family of paths on $(M,d)$ that connects $x$ and $y$. 
We denote $\mathcal C^d_{x,y}$ instead of $\mathcal C_{x,y}$ if there exist more than one metric defined on $M$.
$\hat d$ is the intrinsic (extended) metric induced by $d$, we always have that $d\leq \hat d$. 
We say that the metric $d$ is intrinsic if $\hat d=d$ and we have that a metric $\hat d$ is always intrinsic.

\begin{proposition}
\label{comprimentos iguais}
Let $(M,d)$ be a metric space and $\hat d$ be the intrinsic (extended) metric induced by $d$. 
Then $\ell_{\hat d}(\gamma)=\ell_d (\gamma)$ for every rectifiable curve $\gamma$ in $(M,d)$.
\end{proposition}

If $(M,d)$ is a metric space, $\varepsilon>0$ and $x,y \in M$, then an $\varepsilon$-midpoint of $x$ and $y$ is a point $z\in M$ that satisfies $\vert 2d(x,z)-d(x,y)\vert \leq \varepsilon$ and $\vert 2d(y,z) - d(x,y) \vert \leq \varepsilon$. 

We have the following results about intrinsic metrics. Their proofs can be found in \cite{Burago}.

\begin{proposition}
\label{existenciamidpoint}
If $d$ is an intrinsic metric on $M$ and $\varepsilon> 0$,  then every $x,y \in M$ admit an $\varepsilon$-midpoint. 
\end{proposition}

\begin{proposition}
\label{midpoint}
Let $(M,d)$ be a complete metric space. 
If every $x,y\in M$ admit a $\varepsilon$-midpoint (for every $\varepsilon >0$), then $d$ is intrinsic.
\end{proposition}

$C^0$-Carnot-Carath\'eodory-Finsler-metrics are intrinsic because they come from a length structure (see \cite{Burago}).
A particular case of $C^0$-Carnot-Carath\'eodory-Finsler metrics are the $C^0$-Finsler metrics.
Let $M$ be a differentiable manifold and $TM=\{(x,v);x\in M, v\in T_xM\}$ be its tangent bundle. 
A $C^0$-Finsler metric on $M$ is a continuous function $F:TM \rightarrow \mathbb{R}$ such that $F(x,\cdot)$ is a norm on $T_xM$ for every $x\in M$. 
$F$ induces a metric $d_F$ on $M$ given by
\[
d_F(x,y)=\inf_{\gamma \in \mathcal{S}_{x,y}}\ell_F(\gamma),
\]
where
\[
\ell_F(\gamma)=\int_a^b F(\gamma(t),\gamma^\prime(t))dt,
\]
is the length of $\gamma$ in $(M,F)$ and $\mathcal{S}_{x,y}$ is the family of paths on $M$ which are smooth by parts and connects $x$ and $y$. 
Here the family $\mathcal S_{x,y}$ can be replaced by the family of absolutely continuous paths connecting $x$ and $y$ (see \cite{Berestovskii1}).
If $d:M\times M \rightarrow \mathbb{R}$ is a metric on $M$, then we say that $d$ is $C^0$-Finsler if there exist a $C^0$-Finsler metric $F$ on $M$ such that $d=d_F$.
A differentiable manifold endowed with a $C^0$-Finsler metric is a $C^0$-Finsler manifold.

\begin{remark}
\label{outrofinsler}
There are another usual (in fact more usual) definition of Finsler manifold, where $F$ satisfies other conditions  (see, for instance, \cite{BaoChernShen} and \cite{Deng}): $F$ is smooth on $TM-TM_0$, where $TM_0=\{(p,0)\in TM;p \in M\}$ is the zero section, and  $F(p,\cdot)$ is a Minkowski norm on $T_pM$ for every $p\in M$. 
In order to make this difference clear we put the prefix $C^0$ before Finsler.

$C^0$-Finsler metrics are studied, for instance, in \cite{Berestovskii1}, \cite{Berestovskii2}, \cite{Burago} and \cite{Gribanova}.
\end{remark}

\begin{remark}
\label{relacoes Carnot Caratheodory}
Let $M$ be a connected differentiable manifold and $\mathcal D$ be a completely nonholonomic distribution on $M$.
Let $F$ as in the definition of the $C^0$-Carnot-Carath\'eodory-Finsler metric in the introduction.
Observe that there exist smooths sections $\mathbf g_1, \mathbf g_2$ of inner products in $\mathcal D$ such that $\mathbf g_1 \leq F \leq \mathbf g_2$.
In fact, this is clear locally, and in order to see that this holds globally, just use partition of unity.
This implies that balls in $C^0$-Carnot-Carath\'eodory-Finsler metrics are contained and contains a ball in some Carnot-Carath\'eodory metric of the same distribution. 
Ball-box theorem gives a qualitative behavior of balls in sub-Riemannian manifolds (see \cite{Montgomery}) and this behavior depends only on the distribution $\mathcal D$.
Therefore balls in $C^0$-Carnot-Carath\'eodory-Finsler metrics of the same distribution have the same qualitative behavior.

Another consequence of the ball-box theorem is that a Carnot-Carath\'eodory metric correspondent to a non-trivial completely nonholonomic distribution can't the bounded above by the distance function of a Riemannian metric.
Therefore a $C^0$-Carnot-Caratheodory-Finsler metric correspondent to a non-trivial completely nonholonomic distribution can't be bounded above by a $C^0$-Finsler metric.
\end{remark}

The following theorem is the version of the classical Hopf-Rinow theorem for intrinsic metrics (see \cite{Burago}).

\begin{theorem}[Hopf-Rinow-Cohn-Vossen theorem]
\label{hopf rinow intrinsic}
Let $(M,d)$ be a locally compact metric space endowed with an intrinsic metric.
Then the following assertions are equivalent:
\begin{itemize}
\item $(M,d)$ is complete;
\item $B_d[p,r]$ is compact for every $p \in M$ and $r > 0$;
\item Every geodesic (local minimizer parameterized by arclength) $\gamma:[0,a) \rightarrow M$ can be extended to a continuous path $\bar \gamma:[0,a] \rightarrow M$;
\item There exist a $p \in M$ such that every shortest path parameterized by arclength $\gamma:[0,a) \rightarrow M$ satisfying $\gamma(0)=p$ admits a continuous extension $\bar \gamma:[0,a] \rightarrow M$.
\end{itemize}
\end{theorem}

\begin{lemma}
\label{fechado limitado compacto}
Any closed and bounded subset of a Lie group endowed with a right invariant intrinsic metric is compact. 
In particular, a right invariant intrinsic metric on a Lie group is complete.
\end{lemma}

{\it Proof}

\

A right invariant intrinsic metric $\tilde d$ is a $C^0$-Carnot-Carath\'eodory Finsler metric that comes from a continuous right invariant section of norms $F$ in a right invariant completely nonholonomic distribution $\mathcal D$ (see Theorem \ref{Berestovskii theorem}). 
Observe that $\tilde d \geq d_{\mathbf g}$ for some right invariant Riemannian manifold $\mathbf g$.
In fact, choose a right invariant smooth section of inner products $\bar{\mathbf g}$ on $\mathcal D$ such that $\bar{\mathbf g} \leq F$ and then extend $\bar{\mathbf g}$ to the whole tangent spaces resulting in a right invariant Riemannian metric $\mathbf g$.
The completeness of $(M,\mathbf g)$ and the classical Hopf-Rinow theorem for Riemannian manifolds implies that any closed bounded subset of $(M,\mathbf g)$ is compact.
Therefore any closed bounded subset of $(M,\tilde d)$ is compact because $\tilde d \geq d_{\mathbf g}$.

The last statement is due to Theorem \ref{hopf rinow intrinsic}.$\blacksquare$

\

The rest of this section is about the orbits of a family of vector fields and it is used in Section \ref{existencia sx denso} in order to prove the existence of $X$ such that $S_X$ is dense in $G$ (see \cite{Jurdjevic}).

\begin{definition}
\label{orbita}
Let $\mathcal{F}$ be a family of complete vector fields on a differentiable manifold $M$ and $x\in M$. 
For $X \in \mathcal{F}$, denote by $t\mapsto (\exp tX)(x)$ be the integral curve of $X$ such that $(\exp 0X) (x)=x$.
The orbit $G(x)$ of $\mathcal{F}$ through $x$ is the set of points given by $(\exp t_mX_m)(\exp t_{m-1}X_{m-1}) \ldots (\exp t_1X_1)(x)$, with $m\in \mathbb N$, $t_i \in \mathbb{R}$ and $X_i \in \mathcal{F}$, $i=1,\ldots m$. 
\end{definition}

In the conditions of Definition \ref{orbita}, denote by $Lie(\mathcal{F})$ the Lie algebra generated by $\mathcal{F}$ and let $Lie_x(\mathcal{F})$ be the restriction of $Lie(\mathcal{F})$ to the tangent space $T_xM$.

\begin{theorem}[Hermann-Nagano Theorem]
\label{nagano}
Let $M$ be an analytic manifold and $\mathcal{F}$ be a family of analytic vector fields on $M$. Then
\begin{itemize}
\item each orbit of $\mathcal F$ is an analytic submanifold of $M$;
\item if $N$ is an orbit of $\mathcal{F}$, then the tangent space of $N$ in $x$ is given by $\text{Lie}_x(\mathcal F)$.
\end{itemize}

\end{theorem}

\begin{corollary}
\label{nagano para grupos de Lie}
Let $G$ be a connected Lie group and $V$ be a vector subspace of $\mathfrak{g}$ such that the Lie algebra generated by $V$ is $\mathfrak{g}$.
Let $\{v_1,\ldots,v_k\}$ be a basis of $V$ and $\mathcal{F}$ be the set of left (or right) invariant vector fields with respect to $\{v_1,\ldots,v_k\}$.
Then, for every $x \in G$, there exist $m\in \mathbb{N}$, $t_{i_1},\ldots, t_{i_m} \in \mathbb R$ and $v_{i_1},\ldots, v_{i_m} \in X$ such that $x=\exp(t_{i_m}v_{i_m})\ldots \exp(t_{i_1}v_{i_1})$.
\end{corollary}

{\it Proof}

\

It is enough to observe that $G$ is an analytic manifold, that the right (and left) invariant vector fields with respect to $v_1,\ldots, v_k$ are analytic (see \cite{Helgason}) and that $\text{Lie}_e(\mathcal{F})=\mathfrak g$.
Then the orbit of $e$ is an analytic submanifold of $G$ that contains a neighborhood of $e$ (see Theorem \ref{nagano}) and this neighborhood generates $G$ due to the connectedness of $G$.$\blacksquare$

\section{Convegence of the sequence of induced Hausdorff metrics on Lie groups}
\label{sequencia de metricas}

Let $d$ be a metric on a Lie group $G$ such that $\tau_d=\tau_G$, $\varphi: G \times (G,d) \rightarrow (G,d)$ be the product of $G$ and $X= \{ x_1, \ldots, x_k \} \ni e$ be a finite subset of $G$ such that $H_X=\{e\}$.
In this section, we study properties of the sequence of induced Hausdorff metrics $d,d^1,d^2,\ldots,d^i, \ldots$ on $G$.

For every $j=1, \ldots, k$, we define the metric $d_j(p,q)=d(px_j,qx_j)$ on $G$.
$\tau_{d_j}=\tau_d$ because right translations are homeomorphisms on $G$.
Define also  the metric $d_M(p,q):=\max\limits_{j=1,\ldots, k}d_j(p,q):=\max\limits_{j=1,\ldots, k}d(px_j,qx_j)$ on $G$.

First of all we prove that $\tau_{d_X}=\tau_d$ (Proposition \ref{dedXmesmatopologia}).

\begin{lemma}
\label{comparadXdMd}
$d\leq d_M$ and $d_X \leq d_M$.
\end{lemma}

{\it Proof}

\

The first inequality is obvious. The second inequality follows because
\[
d_X(p,q)
=\max \left\{ \max_i \min_j d (px_i,qx_j), \max_j \min_i d (px_i,qx_j) \right\}
\]
\[
\leq \max \left\{ \max_i d (px_i,qx_i), \max_j d (px_j,qx_j) \right\}
=\max_i d (px_i,qx_i)
=d_M(p,q).\blacksquare
\]

\begin{lemma}
\label{localmenteigual}
For every $g\in G$, there exist a $G$-neighborhood $V$ of $g\in G$ such that $d_X\vert_{V \times V}=d_M\vert_{V \times V}$.
\end{lemma}

{\it Proof}

\

For $i,j=1,\ldots , k$, $i\neq j$, define $\rho_{ij}:G\times G\rightarrow \mathbb R$ as $\rho_{ij}(p,q)=d(px_i,qx_j)-\max_k d(px_k,qx_k)$.
Then for every $g\in G$, there exist a $G$-neighborhood $V$ of $g$ such that 
\[
\rho_{ij}(p,q)>0 \text{ for every } p,q \in V \text{ and every }i\neq j,
\]
because $\rho_{ij}(g,g)>0$.
It implies that
\[
d_X(p,q)
=\max \left\{ \max_i \min_j d (px_i,qx_j), \max_j \min_i d (px_i,qx_j) \right\}
\]
\begin{equation}
\label{vizinhanca pequena}
=\max \left\{ \max_i d (px_i,qx_i), \max_j d (px_j,qx_j) \right\}
=\max_i d(px_i,qx_i)=d_M(p,q)
\end{equation}
for every $p,q \in V$.$\blacksquare$

\begin{remark}
The formula $d_X(p,q)=\max_i d(px_i,qx_i)$ for every $p,q$ in a sufficiently small neighborhood of $g$ will be used several times in this work.
\end{remark}

\begin{lemma}
\label{dMdmesmatopologia} $\tau_d = \tau_{d_M}(=\tau_G)$ on $G$.
\end{lemma}

{\it Proof}

\

The inequality $d\leq d_M$ implies that $\tau_d \subset \tau_{d_M}$.

In order to see that $\tau_{d_M} \subset \tau_d$, observe that an open ball $B_{d_M}(p,r)$ can be written as
\[
B_{d_M}(p,r)=\bigcap_{j=1}^k B_{d_j}(p,r),
\]
what implies that it is an open subset of $\tau_d$. Therefore $\tau_{d_M}\subset \tau_d$.$\blacksquare$

\begin{lemma}
\label{BolasdXcompactas}
$B_{d_X}(p,r)$ is contained in a $G$-compact subset of $G$ for a sufficiently small $r > 0$.
\end{lemma}

{\it Proof}

\

Just observe that
\[
B:=\bigcup_{i=1}^k B_d(px_i,r)\supset B_{d_X}(p,r).
\]
In fact, if $y \not \in B$, then
\[
\begin{array}{ccl}
d_X(p,y)
&
=
&
\max\left\{ \max\limits_i \min\limits_j d (px_i,yx_j), \max\limits_j \min\limits_i d (px_i,yx_j) \right\} \\
&
\geq 
&
\max\left\{ \max\limits_i \min\limits_j d (px_i,yx_j), \min\limits_i d (px_i,y.e) \right\} \\
&
\geq
& 
\max\left\{ \max\limits_i \min\limits_j d (px_i,yx_j), r \right\}
\geq r,
\end{array}
\]
where the second inequality holds because $y \not \in B$. 
Then the result follows because $B$ is contained in a $G$-compact subset of $G$ for a sufficiently small $r$.

\begin{proposition}
\label{dedXmesmatopologia}
$\tau_d = \tau_{d_X}$ on $G$.
\end{proposition}

{\it Proof}

\

It is enough to prove that $\tau_{d_X} = \tau_{d_M}$ due to Lemma \ref{dMdmesmatopologia}.

We know that $d_X \leq d_M$ (Lemma \ref{comparadXdMd}), what implies that $\tau_{d_X} \subset \tau_{d_M}$.

In order to prove that $\tau_{d_M} \subset \tau_{d_X}$, we consider an open ball $B_{d_M}(p,r)$ and we prove that there exist an $\varepsilon > 0$ such that $B_{d_X}(p,\varepsilon)\subset B_{d_M}(p,r)$.
Without loss of generality, we can consider $r > 0$ such that $B_{d_X}(p,r)$ has compact closure (see Lemma \ref{BolasdXcompactas}).
Consider a $G$-neighborhood $V$ of $p$ according to Lemma \ref{localmenteigual}.
We can eventually consider a (further) smaller $r$ in such a way that $B_{d_M}(p,r) \subset V$ (see Lemma \ref{dMdmesmatopologia}). Then $B_{d_X}(p,r)\cap V=B_{d_M}(p,r)$ due to Lemma \ref{localmenteigual}.

If $B_{d_X}(p,r)-V$ is the empty subset, then we have that $B_{d_X}(p,r) = B_{d_M}(p,r)$ and we are done. 
Otherwise we consider
\[
2\varepsilon=\inf\limits_{x \in B_{d_X}(p,r)-V}d_X(p,x).
\]
Observe $d_X$ is $d \times d$-continuous because $\tau_{d_X}\subset \tau_{d}$ and $\varepsilon$ is strictly positive because $p$ is not contained in the $G$-compact subset $\overline{B_{d_X}(p,r)-V} $. 
Therefore
\[
B_{d_X}\left(p,\varepsilon\right)
=B_{d_X}\left(p,\varepsilon \right) \cap V
\subset B_{d_M}(p,r)
\]
what settles the proposition.$\blacksquare$

\begin{definition}
\label{semigrupo}
Let $X=\{x_1,\ldots, x_k\}\ni e$ be a finite subset of a Lie group $G$. The semigroup generated by $X$ is defined as 
\[
S_X=\{x_{i_1}\ldots x_{i_m};m\in \mathbb N, i_j \in (1,\ldots, k), j\in (1,\ldots, m)\}.
\]
\end{definition}

In what follows, if $X=\{x_1,\ldots, x_k\}$, then $X^{-1}:=\{x_1^{-1}, \ldots, x_k^{-1}\}$.

\begin{proposition}
\label{x e menos x}
Let $X$ be a finite subset of a Lie group $G$. 
Then $S_X$ is dense in $G$ iff $S_{X^{-1}}$ is dense in $G$.
\end{proposition}

{\it Proof}

\

It is enough to observe that if $i:G \rightarrow G$ is the inversion map, then $i(S_X)=S_{X^{-1}}$.$\blacksquare$

\

\begin{lemma}
\label{d1 e maior}
Let $G$ be a Lie group and
$d:G\times G\rightarrow \mathbb R$ be a metric on $G$ such that $\tau_d=\tau_G$. 
Let $\varphi: G \times (G,d)\rightarrow (G,d)$ be the product of $G$ and
$X=\{x_1,\ldots, x_k\}\ni e$ be a finite subset of $G$ such that $H_X=\{e\}$.
Then $d\leq d^1$ and the sequence of induced Hausdorff metrics is increasing.
\end{lemma}

{\it Proof}

\

Let $x,y \in G$ and $\gamma:[a,b]\rightarrow G$ be a $d_X$-path (which is also a $G$-path due to Proposition \ref{dedXmesmatopologia}) connecting $x$ and $y$ (If there isn't any path connecting $x$ and $y$, then $d(x,y)\leq d^1(x,y) = \infty$).
Then
\[
\ell_{d_X}(\gamma)=\sup_{\mathcal P}\sum_{i=1}^{n_{\mathcal P}}d_X(\gamma(t_i),\gamma(t_{i-1})).
\]
Cover $\gamma([a,b])$ with open subsets $V$ of Lemma \ref{localmenteigual} such that (\ref{vizinhanca pequena}) holds. 
Let $\varepsilon$ be the $d$-Lebesgue number of this covering. Let $\delta >0$ such that if $\vert t_1 - t_2\vert < \delta$, then $d(\gamma(t_2),\gamma(t_1))<\varepsilon$.
Consider only partitions with norm less that $\delta$. Then
\begin{equation}
\label{calculo comprimento}
\ell_{d_X}(\gamma)
=
\sup_{\mathcal P}\sum_{i=1}^{n_{\mathcal P}}\max_j d(\gamma(t_i)x_j,\gamma(t_{i-1})x_j) 
\geq
\sup_{\mathcal P}\sum_{i=1}^{n_{\mathcal P}} d(\gamma(t_i),\gamma(t_{i-1})) = \ell_d(\gamma)
\end{equation}
due to (\ref{vizinhanca pequena}).
Then
\[
d^1(x,y)
=\inf_{\gamma \in \mathcal C_{x,y}^{d_X}} \ell_{d_X}(\gamma)
\geq \inf_{\gamma \in \mathcal C_{x,y}^d} \ell_{d}(\gamma)
=\hat d(x,y)\geq d(x,y).\blacksquare
\]

\begin{theorem}
\label{local maximo 1}
Let $G$ be a Lie group and
$d:G\times G\rightarrow \mathbb R$ be a metric on $G$ such that $\tau_d=\tau_G$.
Let $\varphi: G \times (G,d)\rightarrow (G,d)$ be the product of $G$ and
$X=\{x_1,\ldots, x_k\}\ni e$ be a finite subset of $G$ such that $H_X=\{e\}$.
Suppose that $d$ is bounded above by a right invariant intrinsic metric $\mathbf d$ (what implies that $G$ is connected). 
Then every $d^i$ is bounded above by $\mathbf d$ and the sequence of induced Hausdorff metrics converges pointwise to a metric $d^\infty$.
\end{theorem}

{\it Proof}

\

We prove that if $d\leq \mathbf d$, then $d^1 \leq \mathbf d$.
This is enough to prove that $d^i\leq \mathbf d$ for every $i$: just  iterate the process, replacing $d$ by $d^{i-1}$.

Let $\gamma$ be a path defined on a closed interval. Then
\begin{eqnarray}
\ell_{d_X}(\gamma)
&
=
&
\sup_{\mathcal P}\sum_{i=1}^{n_{\mathcal P}}\max_j d(\gamma(t_i)x_j,\gamma(t_{i-1})x_j) \nonumber \\
&
\leq
& 
\sup\limits_{\mathcal P}\sum\limits_{i=1}^{n_{\mathcal P}}\max\limits_j \mathbf d(\gamma(t_i)x_j,\gamma(t_{i-1})x_j) \nonumber
\\
&
=
&
\sup\limits_{\mathcal P}\sum\limits_{i=1}^{n_{\mathcal P}} \mathbf d (\gamma(t_i),\gamma(t_{i-1}))= \ell_{\mathbf d}(\gamma) \label{majoracao}
\end{eqnarray}
because $\mathbf d$ is right invariant. 
From (\ref{majoracao}) and the fact that $\mathbf d$ is intrinsic, we have that 
\[
d^1(x,y) 
= 
\inf_{\gamma \in \mathcal{C}_{x,y}^{d_X}} \ell_{d_X}(\gamma) 
\leq 
\inf_{\gamma \in \mathcal{C}_{x,y}^{\mathbf d}} \ell_{\mathbf d}(\gamma)
=
\mathbf d(x,y)
\]
for every $x,y \in G$ (Observe that $\mathcal C_{x,y}^{d_X}= \mathcal C_{x,y}^{\mathbf d}$ due to Proposition \ref{dedXmesmatopologia}).

Observe that $d^i$ is an increasing sequence of metrics bounded above by $\mathbf d$.
Then $d^i$ converges pointwise to a function $d^\infty$ and it is easy to prove that $d^\infty$ is a metric.$\blacksquare$

\

Observe that if $d$ is a metric on a group $G$, then $\bar d: G \times G \rightarrow \mathbb R$, defined as $\bar d(x,y)=\sup_{\sigma \in G}d(x\sigma,y\sigma)$, is a right invariant extended metric. 
Its intrinsic extended metric $\hat{\bar d}(x,y)=\inf_{{\mathcal C}^{\bar{d}}_{x,y}} \ell_{\bar d}(\gamma)$ is also a right invariant.

\begin{corollary}
\label{d barra chapeu limitante} 
Let $G$ be a Lie group and
$d:G\times G\rightarrow \mathbb R$ be a metric on $G$ such that $\tau_d=\tau_G$.
Let $\varphi: G \times (G,d)\rightarrow (G,d)$ be the product of $G$ and
$X=\{x_1,\ldots, x_k\}\ni e$ be a finite subset of $G$ such that $H_X=\{e\}$. 
Suppose that $\hat{\bar{d}}$ is a metric on $G$.
Then the sequence of induced Hausdorff metrics converges to a metric $d^\infty$ on $G$.
\end{corollary}

{\it Proof}

\

It is enough to see that $\hat{\bar{d}}$ is a right invariant intrinsic metric such that $d\leq \hat{\bar{d}}$.$\blacksquare$  

\begin{lemma}
\label{translacao a direita decresce}
Let $d^i$ be a sequence of induced Hausdorff metrics converging to $d^\infty$.
Then $d^\infty(x,y) \geq d^\infty(x \sigma,y \sigma)$ for every $\sigma \in \bar S_X$ and every $x,y \in G$.
In particular $\ell_{d^\infty}(\gamma) \geq \ell_{d^\infty}(\gamma \sigma)$ for every $\sigma \in \bar S_X$ and every path $\gamma$. 
Moreover if $\bar S_X$ is a subgroup, then $d^\infty$ is invariant by the right action of $\bar S_X$. 
In particular, if $\bar S_X=G$, then $d^\infty$ is right invariant.
\end{lemma}

{\it Proof}

\

In order to prove that $d^\infty(x,y) \geq d^\infty(x\sigma,y \sigma)$ for every $(x,y,\sigma) \in G \times G \times \bar S_X$, it is enough to prove that $d^\infty(x,y) \geq d^\infty(x x_j,y x_j)$ for every $(x,y,x_j) \in G \times G \times X$ because every $\sigma \in \bar S_X$ can be arbitrarily approximated by product of elements of $X$.

Let $d$ be a metric on $G$.
Then 

\begin{eqnarray}  
\ell_{d^1}(\gamma) &
= & \sup\limits_{\mathcal P}\sum\limits_{i=1}^{n_{\mathcal P}} d^1(\gamma(t_i),\gamma(t_{i-1})) 
= 
\sup\limits_{\mathcal P}\sum\limits_{i=1}^{n_{\mathcal P}} d_X(\gamma(t_i),\gamma(t_{i-1})) \nonumber \\
&
=
&
\sup\limits_{\mathcal P}\sum\limits_{i=1}^{n_{\mathcal P}} \max\limits_j d(\gamma(t_i)x_j,\gamma(t_{i-1})x_j) \nonumber \\
& 
\geq
&
\sup\limits_{\mathcal P}\sum\limits_{i=1}^{n_{\mathcal P}} d(\gamma(t_i)x_j,\gamma(t_{i-1})x_j)=\ell_{d}(\gamma x_j) \label{multiplica xj diminui}
\end{eqnarray}
for every $x_j \in X$.
The second equality holds due to Proposition \ref{comprimentos iguais}.

Fix $(x,y,x_j) \in G \times G \times X$. Observe that $d^i(x,y)\geq d^{i-1}(x x_j, y x_j)$ for $i \geq 2$ due to (\ref{multiplica xj diminui}), $\lim_{i \rightarrow \infty}d^i(x,y)=d^\infty(x,y)$ and $\lim_{i\rightarrow \infty} d^{i-1}(x x_j, y x_j)=d^\infty (xx_j, yx_j)$.
Then $d^\infty(x,y)\geq d^\infty(xx_j, yx_j)$, what implies that $d^\infty(x,y)\geq d^\infty(x\sigma, y\sigma)$ for every $\sigma \in \bar S_X$.

In order to prove that if $\bar S_X$ is a subgroup, then $d^\infty$ is invariant by the right action of $\bar S_X$, it is enough to observe that $d^\infty(x,y)\geq d^\infty(x\sigma, y\sigma)\geq d^\infty(x\sigma \sigma^{-1}, y\sigma\sigma^{-1})=d^\infty(x,y)$.$\blacksquare$

\

Now we prove that if $d$ is complete and bounded above by a right invariant intrinsic metric, then $d^\infty$ is intrinsic (Theorem \ref{completo eh intrinseco}). 

\begin{lemma}
\label{osdoiscompletos}
Let $M$ be a set and suppose that $d$ and $\rho$ are metrics on $M$ such that $d \leq \rho$ and $\tau_d=\tau_\rho$. 
If $d$ is complete, then $\rho $ is also complete. 
\end{lemma}

{\it Proof}

\

Just notice that every Cauchy sequence in $(M,\rho)$ is a Cauchy sequence in $(M,d)$. 
Then the sequence converges in both metrics spaces because $\tau_d=\tau_\rho$.$\blacksquare$

\begin{lemma}
\label{limiteintrinseco}
Let $M$ be locally compact topological space and suppose that $d^i$ is a sequence intrinsic metrics on $M$ that converges pointwise to metric  $d^\infty$.
Suppose that there exist a complete intrinsic metric $d_l$ and a metric $d_h$ on $M$ such that $\tau_{d_l}=\tau_{d_h}=\tau_M$ and $d_l \leq d^i \leq d_h$ for every $i$.
Then $d^\infty$ is a complete intrinsic metric.
\end{lemma}

{\it Proof}

\

First of all $d^\infty$ is complete due to Lemma \ref{osdoiscompletos} and because $d^\infty \geq d_l$. Moreover it is straightforward that $\tau_{d^\infty}=\tau_M$.

Let $\varepsilon >0$ and $x,y \in M$. 
We will prove that $x$ and $y$ admits an $\varepsilon$-midpoint $z$ with respect to the metric $d^\infty$ (see Proposition \ref{midpoint}).

Let $z_i$ be an $\varepsilon /3$-midpoint of $x$ and $y$  with respect to the metric $d_i$ (see Proposition \ref{existenciamidpoint}). We claim that the sequence $z_i$ is contained in a compact subset.
In fact 
\[
\left\vert 2d^i(x,z_i)-d^i(x,y) \right\vert \leq \frac{\varepsilon}{3}
\]
implies that
\[
2d^i(x,z_i) 
\leq d^i(x,y) + \frac{\varepsilon}{3}
\]
and
\[
2d_l(x,z_i) 
\leq 2d^i(x,z_i) 
\leq d^i(x,y) + \frac{\varepsilon}{3} 
\leq d_h(x,y) + \frac{\varepsilon}{3}.
\]
Then the sequence $z_i$ is bounded in $(M,d_l)$ and it is contained in a compact subset due to Theorem \ref{hopf rinow intrinsic}.

Let $z$ be an accumulation point of $z_i$.
We claim that $z$ is an $\varepsilon$-midpoint with respect to $d^\infty$.
First of all, taking a subsequence $z_i$ coverging to $z$, we claim that $d^i(x,z_i)$ converges to $d^\infty (x,z)$ as $i$ goes to infinity.
In fact, for every $\varepsilon >0$, there exist $N_1\in \mathbb N$ such that
\[
\vert d^\infty(x,z) - d^i(x,z_i) \vert 
\leq \vert d^\infty(x,z) - d^i(x,z) + d^i(x,z)- d^i(x,z_i) \vert
\]
\[
\leq \vert d^\infty(x,z) - d^i(x,z)\vert + \vert d_h(x,z)- d_h(x,z_i) \vert < \varepsilon /3
\]
for every $i\geq N_1$, and $d^i(x,z_i)$ converges to $d^\infty(x,z)$. Then for every $\varepsilon >0$, there exist a $N_2 \in \mathbb{N}$ such that
\[
\vert 2d^\infty(x,z) - d^\infty(x,y)\vert
\]
\[
\leq \left\vert 2d^\infty(x,z)-2d^i(x,z_i) \vert+ \vert 2d^i(x,z_i)-d^i(x,y) \vert + \vert d^i(x,y)-d^\infty(x,y) \right\vert 
< \varepsilon
\]
for every $i \geq N_2$. The inequality 
\[
\vert 2d^\infty(y,z) - d^\infty(x,y)\vert < \varepsilon
\]
is analogous. Then $z$ is an $\varepsilon$-midpoint of $x$ and $y$ in $(M,d^\infty)$ and $d^\infty$ is intrinsic due to Proposition \ref{midpoint}.$\blacksquare$

\begin{theorem}
\label{completo eh intrinseco}
Let $G$ be a Lie group and
$d:G\times G\rightarrow \mathbb R$ be a complete metric on $G$ such that $\tau_d=\tau_G$.
Let $\varphi: G \times (G,d)\rightarrow (G,d)$ be the product of $G$ and
$X=\{x_1,\ldots, x_k\}\ni e$ be a finite subset of $G$ such that $H_X=\{e\}$. 
Suppose that there exist an intrinsic right invariant metric $\mathbf d$ on $G$ such that $d \leq \mathbf d$.
Then the sequence of induced Hausdorff metrics converges to a complete and intrinsic metric $d^\infty$. 
\end{theorem}

{\it Proof}

\

Observe that $\tau_{\mathbf d} = \tau_d = \tau_G$ because $\mathbf d$ is $C^0$-Carnot-Carath\'eodory-Finsler and $d\leq d^i \leq \mathbf d$.
Notice that the sequence $(d^i)_{i \in \mathbb N}$ converges to a metric $d^\infty$ due to Theorem \ref{local maximo 1} and that we are in the conditions of Lemma \ref{limiteintrinseco}, with $d_h = \mathbf d$ and $d_l=d^1$.
Therefore $d^\infty$ is complete and intrinsic.$\blacksquare$

\section{Existence of dense semigroups $S_X$}
\label{existencia sx denso}

Let $G$ be a connected Lie group. 
In this section we prove the existence of a finite subset $X=\{x_1,\ldots, x_k\}\ni e$ of $G$ such that $H_X=\{e\}$ and $\bar S_X = G$.

We begin with an example.

\begin{example}
\label{denso na reta}
Let $G=\mathbb R$ be the additive group of real numbers and consider $X=\{-1,0,\sqrt 2\}$.
We claim that $S_X$ is dense.
Let $q \in \mathbb{R}$ and $\varepsilon >0$.
We will find a $q_\varepsilon \in S_X$ such that $\vert q - q_\varepsilon \vert < \varepsilon$.

It is easy to see that there are a sequence $p_1,p_2,\ldots$ such that $p_{i+1}-p_i = -1$ or $p_{i+1}-p_i=\sqrt{2}$ and a infinite number of ${p_i}'s$ are in $[0,1]$.
Therefore there exist $p_j$ and $p_k$ in $S_X$, $j<k$, such that $\vert p_k-p_j\vert<\varepsilon$.
Set $p_\varepsilon := p_k-p_j \in S_X$.
In order to fix ideas, suppose that $p_\varepsilon >0$. 
Then there exist an element $p \in S_X$ such that $p<q$ and a $l\in \mathbb N$ such that $q_\varepsilon:=p+l.p_\varepsilon$ satisfies $\vert q- q_\varepsilon \vert < \varepsilon$. The case $p_\varepsilon <0$ is analogous.
Then $S_X$ is dense in $\mathbb R$.
We also have that $S_{X^{-1}}$ is dense in $G$ due to Proposition \ref{x e menos x}.

This example is easily generalized to the case $G=\mathbb R^n$ and $X=\{-1,0,\sqrt 2\}\times \ldots \times \{-1,0,\sqrt 2\}$. $\blacksquare$

\end{example}

Example \ref{denso na reta} can be generalized to Lie groups in the following way:

\begin{theorem}
\label{grupos de Lie}
Let $G$ be a connected Lie group and $\mathfrak g$ be its Lie algebra.
Let $V$ be a vector subspace of $\mathfrak g$ such that the Lie algebra generated by $V$ is $\mathfrak{g}$. Let $\{v_1,\ldots, v_k\}$ be a basis of $V$ such that 
\[
\left\{ \sum_{i=1}^k a_iv_i, a_i \in [-2,2], i=1,\ldots, k \right\}
\]
is contained in a neighborhood $A$ of $0 \in \mathfrak{g}$ such that $\exp\vert_A$ is a diffeomorphism over its image.
Consider 
\[
\tilde X=\{-v_1,-v_2,\ldots,-v_k, 0 , \sqrt{2}v_1, \ldots, \sqrt{2}v_k\}.
\]
Denote $x_i=\exp(-v_i)$, $y_i=\exp(\sqrt{2}v_i)$ and $X=\{e,x_1,\ldots,x_k,y_1,\ldots, y_k\}$.
Then $\bar S_X=G$.
\end{theorem} 

{\it Proof}

\

Let $g \in H_X$. 
Then $gX=X$ and $e \in X$ implies that $g \in X$ and $g^{-1} \in X$. 
If $g \neq e$, then $g=\exp(-v_i)$ or $g=\exp (\sqrt{2}v_i)$.
In order to fix ideas, suppose that $g=\exp(-v_i)$. 
Then we have that $g^{-1}=\exp(v_i)\in X$, which isn't possible because $\exp\vert_A$ is a diffeomorphism over its image.
The other case is analogous.
Consequently $H_X=\{e\}$.

Let us prove that $S_X$ is dense in $G$. Example \ref{denso na reta} states that $S_{\tilde X} \cap \mathbb{R}v_i$ is dense in $\mathbb{R}v_i$ for every $i=1,\ldots,n$.
This implies that $S_X \cap \exp(\mathbb Rv_i)$ is dense in $\exp(\mathbb{R}v_i)$.
By Corollary \ref{nagano para grupos de Lie}, for every $x\in G$, there exist $m \in \mathbb{N}$, $t_{i_1},\ldots, t_{i_m} \in \mathbb R$ and $v_{i_1},\ldots, v_{i_m} \in X$ such that $x=\exp(t_{i_m}v_{i_m})\ldots \exp(t_{i_1}v_{i_1})$.
But we know that for every $j=1,\ldots, m$, there exist a sequence of points in $S_X$ converging to $\exp(t_{i_j}v_{i_j})$.
Therefore we have a sequence of points in $S_X$ converging to $x$ and $S_X$ is dense in $G$.$\blacksquare$

\section{The case $\bar S_X=G$}
\label{Finlser fomula}

In this section we study the case where $\bar S_X=G$ and $d$ is complete. 
We prove that the sequence of induced Hausdorff metrics converges to a $C^0$-Carnot-Carath\'eodory-Finsler metric $d^\infty$. 
Moreover we give a necessary and sufficient condition in order to $d^\infty$ be $C^0$-Finsler. 
In this case, an explicit formula for the Finsler metric is obtained. 
All these results are in Theorem \ref{teorema principal}.

\begin{lemma}
\label{d barra chapeu e fracao}
Let $G$ be a connected Lie group endowed with a metric $d$ such that $\tau_d=\tau_G$. Suppose that $\hat{\bar{d}}(x,y)= \infty$ for some $x,y \in G$. 
Then there exist a $v \in \mathfrak g$ such that 
\[
\sup_{\sigma \in G}\sup_{t \neq 0} \frac{d(\exp(tv) \sigma, \sigma)}{\vert t\vert}
 =\infty.
\]
\end{lemma}

{\it Proof}

\

If there exist a $\delta > 0$ such that $\hat{\bar d}(\exp(tv),e)$ is finite for every $v \in \mathfrak g$ and $t \in (-\delta,\delta)$, then it is straightforward that $\hat{\bar{d}}(x,y)$ is finite for every $x,y \in G$ due to the right invariance of $\hat{\bar d}$ and the connectedness of $G$. 
Then if $\hat{\bar{d}}(x,y)= \infty$ for some $x,y \in G$, then there exist $v \in \mathfrak g$ and $t\in \mathbb R$ such that $\hat{\bar d}(\exp(tv),e) = \infty$.
It implies that for every $C>0$, there exist a partition $\mathcal P$ of $[0,t]$ such that
\[
\sum_{i=1}^{n_{\mathcal P}}\bar d(\exp(t_iv),\exp(t_{i-1}v)) 
= 
\sum_{i=1}^{n_{\mathcal P}}\bar d(\exp(t_i-t_{i-1})v,e) > Ct
\]  
due to the right invariance of $\bar d$. Then the average of
\[
\frac{\bar d(\exp(t_i-t_{i-1})v,e)}{t_i - t_{i-1}}
\]
is greater than or equal to $C$. Therefore there exist $t > 0$ such that
\[
\frac{\bar d(\exp(tv,e)}{t}>C.
\]
Thus 
\[
\sup_{\sigma \in G}\sup_{t \neq 0} \frac{d(\exp(tv) \sigma, \sigma)}{\vert t\vert} = \sup_{t \neq 0}\frac{\bar d(\exp(tv),e)}{t}
 =\infty. \blacksquare
\]

\begin{lemma}
\label{supelimsup}
Let $G$ be a Lie group, $d$ be a metric on $G$ satisfying $\tau_d = \tau_G$ and $v \in \mathfrak g$.
Then
\begin{equation}
\label{igualdade limsup sup 2}
\sup_{\sigma \in G}\limsup_{t \rightarrow 0} \frac{d(\exp(tv) \sigma, \sigma)}{\vert t\vert}
 =\sup_{\sigma \in G}\sup_{t \neq 0} \frac{d(\exp(tv) \sigma, \sigma)}{\vert t\vert}.
\end{equation}

\end{lemma}

{\it Proof}

\

The inequality $\leq $ follows directly from the definition of $\limsup$. 
Let us prove the inequality $\geq$.
If the left hand side of (\ref{igualdade limsup sup 2}) is infinity, then there is nothing to prove.
Suppose that it is equal to $L \in \mathbb R$.
It implies that 
\[
\limsup_{t \rightarrow 0} \frac{d(\exp(tv)\sigma, \sigma)}{\vert t\vert}\leq L
\]
for every $\sigma \in G$.
Then for every $\varepsilon >0$ and $\sigma \in G$, there exist a $\delta >0$ such that 
\[
d(\exp(tv)\sigma, \sigma)
\leq (L + \varepsilon)\vert t \vert
\]
whenever $t \in (-\delta, \delta)$. 

Now fix $\sigma \in G$ and $t > 0$.
Define $\gamma:[0,t] \rightarrow G$ as $\gamma(s) = \exp(sv).\sigma$.
We can find a partition $\mathcal P=\{0=t_0 < t_1 < \ldots < t_k = t\}$ such that 
\begin{equation}
\label{limita inclinacao}
d(\exp(t_{i}v)\sigma,\exp(t_{i-1}v)\sigma) \leq (L + \varepsilon)(t_i - t_{i-1}).
\end{equation}
In fact, for every $s\in [0,t]$, we choose $\delta_s >0 $ such that  
\[
d(\exp(\tau v)\sigma,\exp(sv)\sigma) \leq (L + \varepsilon)\vert \tau - s \vert.
\] 
for every $\tau \in I_s:=[0,t] \cap (s -\delta_s, s + \delta_s)$. 
If $s \in (0,t)$, we can choose $I_s = (s -\delta_s, s + \delta_s)$.
From $\{I_s\}_{s \in [0,t]}$, we can choose a finite subcover $\{I_{\tilde k}\}$. 
From this finite subcover, we drop a $I_{\tilde k_1}$, once a time, if there exist $I_{\tilde k_2}$ such that $I_{\tilde k_1} \subset I_{\tilde k_2}$.
We end up with a family such that one element is not contained in another element. 
Denote this family by $\{[0,0 + \delta_0),(t_2 - \delta_{t_2},t_2 + \delta_{t_2}),(t_4 - \delta_{t_4}, t_4 + \delta_{t_4}), \ldots,(t-\delta_t,t]\}$, with $t_{i-2} < t_i$ for every $i$ ($t_{i-2} = t_i$ doesn't happen).
It is not difficult to see that $0 < t_2-\delta_{t_2} < t_4 - \delta_{t_4} < \ldots < t - \delta_t$  and $0 + \delta_0 < t_2 + \delta_{t_2} < t_4 + \delta_{t_4} < \ldots < t$. 
Now we choose the ``odd'' points. 
$t_1 \in (0,t_2)$ is chosen in $[0,0 + \delta_0) \cap (t_2 - \delta_{t_2},t_2 + \delta_{t_2})$. 
$t_3 \in (t_2,t_4)$ is chosen in $ (t_2 - \delta_{t_2},t_2 + \delta_{t_2}) \cap (t_4 - \delta_{t_4}, t_4 + \delta_{t_4}) $ and so on.
Then $\mathcal P=\{0=t_0 < t_1 < \ldots < t_k = t\}$ is such that (\ref{limita inclinacao}) is satisfied.

Therefore 
\[
\begin{array}{ccl}
d(\exp(tv)\sigma,\sigma) & \leq & \sum\limits_{i=1}^{n_\mathcal P} d(\exp(t_{i}v)\sigma,\exp(t_{i-1}v)\sigma) \leq (L+\varepsilon)\vert t\vert
\end{array}
\]
for every $\sigma$ and $t>0$ and we can conclude that
\[
\sup_{\sigma \in G}\sup_{t > 0} \frac{d(\exp(tv) \sigma, \sigma)}{\vert t\vert} \leq L.
\]

If $t<0$, then
\begin{eqnarray}
\sup_{\sigma \in G}\sup_{t < 0} \frac{d(\exp(tv) \sigma, \sigma)}{\vert t\vert} & = & \sup_{\sigma \in G}\sup_{t > 0} \frac{d(\exp(t(-v)) \sigma, \sigma)}{\vert t\vert} \nonumber \\ 
& \leq & \sup_{\sigma \in G}\limsup_{t \rightarrow 0} \frac{d(\exp(t(-v)) \sigma, \sigma)}{\vert t\vert} \nonumber \\
& = & \sup_{\sigma \in G}\limsup_{t \rightarrow 0} \frac{d(\exp(tv) \sigma, \sigma)}{\vert t\vert} = L,
\end{eqnarray}
what settles the proposition.$\blacksquare$

\

In what follows we identify $G \times \mathfrak g$ with $TG$ through the correspondence $(g,v) \mapsto (g,dR_g(v))$, where $R_g$ denotes the right translation by $g$ on $G$.

\begin{theorem}
\label{teorema principal}
Let $G$ be a connected Lie group endowed with a complete metric $d$ such that $\tau_d=\tau_G$, $\varphi: G \times G \rightarrow G$ be the product of $G$ and $X=\{x_1,\ldots x_k\} \ni e$ be a finite subset of $G$ such that $H_X=\{e\}$ and $\bar S_X=G$.
\begin{enumerate}
\item If $d$ is bounded above by a right invariant intrinsic metric $\mathbf d$, then the sequence of induced Hausdorff metrics converges to a right invariant $C^0$-Carnot-Carath\'eodory-Finsler metric $d^\infty$.
In particular, if $\hat{\bar{d}}$ is a metric, then $d^\infty$ is a right invariant $C^0$-Carnot-Carath\'eodory-Finsler.
\item If $d$ is bounded above by a right invariant $C^0$-Finsler metric $\mathbf d$, then the sequence of induced Hausdorff metrics converges to a right invariant $C^0$-Finsler metric $d^\infty$.
\item Suppose that the sequence of induced Hausdorff metrics converges pointwise to a metric $d^\infty$.
Then $d^\infty$ is $C^0$-Finsler iff 
\begin{equation}
\label{f tilde}
\tilde F(g,v):=\sup_{\sigma\in G}\limsup_{t \rightarrow 0}\frac{d(\exp(tv)g\sigma,g\sigma)}{\vert t \vert}
\end{equation}
is finite for every $(g,v) \in TG$, and in this case $\tilde F$ is the Finsler metric on $G$.
\end{enumerate}
\end{theorem}

{\it Proof}

\

(1) We know that the sequence of induced Hausdorff metrics converges to a intrinsic metric $d^\infty$ due to Theorem \ref{completo eh intrinseco}. 
It is a right invariant due to Lemma \ref{translacao a direita decresce}.
Therefore $d^\infty$ is $C^0$-Carnot-Carath\'eodory-Finsler due to Theorem \ref{Berestovskii theorem}.

If $\hat{\bar{d}}$ is finite, then $\mathbf d=\hat{\bar{d}}$ is a right invariant intrinsic bound for $d$ and the result follows.

\

(2) This is a particular case of item (1) and it is direct consequence of the second paragraph of Remark \ref{relacoes Carnot Caratheodory}.

\

(3) First of all we prove that if $\tilde F$ is finite, then $d^\infty$ is $C^0$-Finsler.

If $\tilde F$ is finite, then $\hat{\bar d}$ is also finite due to Lemmas \ref{d barra chapeu e fracao} and \ref{supelimsup}.
Thus $d^\infty$ is a right invariant $C^0$-Carnot-Carath\'eodory-Finsler due to item (1).

It is a direct consequence of the ball-box theorem that if $d^\infty$ is a Carnot-Carath\'eodory metric on a differentiable manifold $M$ with respect to a proper smooth distribution $\mathcal D$, then for every $p \in M$ and $v \not\in \mathcal D_p$ we have that
\[
\lim_{t \rightarrow 0}\frac{d^\infty (\gamma(t), p)}{\vert t \vert}=\infty,
\]
where $\gamma:(-\varepsilon,\varepsilon) \rightarrow M$ is a smooth path such that $\gamma(0)=p$ and $\gamma^\prime(0)=v$.
The same conclusion holds for $C^0$-Carnot-Carath\'eodory-Finsler metrics due to Remark \ref{relacoes Carnot Caratheodory}.
Then $d^\infty$ is $C^0$-Finsler whenever $\tilde F$ is finite.

Now we prove that if $d^\infty$ is $C^0$-Finsler, then $\tilde F$ is finite.

If $F:TG \rightarrow \mathbb R$ is the $C^0$-Finsler metric correspondent to $d^\infty$, then
\[
F(g,v)
=\lim_{t\rightarrow 0} \frac{d^\infty(\exp(tv)g,g)}{\vert t\vert} 
= \sup_{\sigma \in G}\lim_{t\rightarrow 0} \frac{d^\infty(\exp(tv)g \sigma,g \sigma)}{\vert t\vert} \geq \tilde F(g,v),
\]
where the first equality is due to Theorem 3.7 of \cite{BenettiFukuoka} and the second equality is due to the right invariance of $d^\infty$.
Therefore $\tilde F$ is finite.
Observe that this part of the proof also shows that $\tilde F \leq F$ if $d^\infty$ is $C^0$-Finsler.

Finally we show that if $d^\infty$ is $C^0$-Finsler, then $F \leq \tilde F$.

Fix $t>0$ and $v \in \mathfrak g$. Define $\gamma:[0,t] \rightarrow G$ by $\gamma(s)=\exp(sv)$ and notice that if $\xi \in G$, then

\begin{eqnarray}
d^1(\exp(tv)\xi, \xi) 
& 
\leq 
&
\ell_{d_X} (\gamma \xi)
=  
\sup\limits_{\mathcal P} \sum\limits_{i=1}^{n_{\mathcal P}} d_X(\exp(t_i v)\xi,\exp(t_{i-1}v)\xi) \nonumber \\
& 
= 
& \sup\limits_{\mathcal P} \sum\limits_{i=1}^{n_{\mathcal P}} \max\limits_j d(\exp(t_i v)\xi x_j,\exp(t_{i-1}v)\xi  x_j) \nonumber \\
&
= 
& \sup\limits_{\mathcal P} \sum\limits_{i=1}^{n_{\mathcal P}} \max\limits_j d(\exp((t_i-t_{i-1}) v)h_j,h_j) \nonumber \\
&
\leq
& \sup\limits_{\mathcal P} \sum\limits_{i=1}^{n_{\mathcal P}} \sup\limits_{\sigma \in G} d(\exp((t_i-t_{i-1}) v) \sigma, \sigma) \nonumber
\end{eqnarray}

where $h_j = \exp(t_{i-1}v)\xi x_j$. Then
\[
\begin{array}{ccl}
d^1(\exp(tv)\xi, \xi) 
& 
\leq 
&
\sup\limits_{\mathcal P} \sum\limits_{i=1}^{n_{\mathcal P}} \tilde F(v)(t_i-t_{i-1}) = \tilde F(v)t,
\end{array}
\]
where the inequality is due to Lemma \ref{supelimsup}. 
If we iterate this process, replacing $d^1$ by $d^2, d^3, \ldots$, we get 
\begin{equation}
\label{compara com is}
d^i(\exp(tv)\xi,\xi) \leq \tilde F(v) t
\end{equation} 
for every $i \in \mathbb N$, $t > 0$ and $\xi \in G$. 
But $\tilde F(v)=\tilde F(-v)$ because
\begin{eqnarray}
\tilde F(v)
&
=
&
\sup_{\sigma\in G}\limsup_{t \rightarrow 0}\frac{d(\exp(tv)\sigma,\sigma)}{\vert t \vert} \nonumber \\
&
=
&
\sup_{\sigma\in G}\limsup_{t \rightarrow 0}\frac{d(\exp(tv)\exp(-tv)\sigma,\exp(-tv)\sigma)}{\vert t \vert} 
= \tilde F(-v),
\end{eqnarray}
and
\begin{equation}
\label{d infinito e f til}
d^\infty(\exp(tv)\sigma,\sigma) \leq \tilde F(v)\vert t \vert
\end{equation}
holds for every $t \in \mathbb R$ and $\sigma \in G$ due to (\ref{compara com is}). Therefore
\[
F(v)
=
\lim_{t \rightarrow 0}\frac{d^\infty(\exp(tv),e)}{\vert t \vert}\leq \tilde F(v),
\]
what settles the theorem.$\blacksquare$

\section{Further examples}
\label{secao exemplos}

\begin{example}
Consider the additive group $G=\mathbb R$ with the finite metric $d(x,y)=\vert\arctan (x) - \arctan (y) \vert$.
Let $X=\{ -1,0,\sqrt 2\}$.
Then $\bar S_X=\mathbb R$ (see Example \ref{denso na reta}).
Its maximum derivative is equal to one and it is not difficult to see that $\hat{\bar d}$ is a metric on $G$.
Then $d^\infty$ is an invariant $C^0$-Finsler metric, with $F(\cdot)$ equal to the Euclidean norm (See Theorem \ref{teorema principal}).
Therefore $d^\infty$ is the Euclidean metric.

Now we prove that every $d^i$ of the sequence of induced Hausdorff metrics is a finite metric.
First of all notice that
\[
\frac{d}{dt}\arctan t = \frac{1}{1+t^2}.
\]
It means that if $0<x<y$, then 
\[
d(x,y) = \arctan y - \arctan x \leq \frac{1}{1+x^2} (y - x).
\]
because $t \mapsto \arctan t$ has the concavity pointed downwards for $t >0$.
For $1<x<y$, we have that
\[
d^1(x,y)=\sup_{\mathcal P}\sum_{i=1}^{n_{\mathcal P}} \max_j d(x+x_j,y+x_j) \leq \frac{1}{1+(x-1)^2}(y-x).
\]
We can iterate this procedure and show that if $i \in \mathbb N$ and $i < x< y$, then 
\[
d^i(x,y) \leq \frac{1}{1+(x-i)^2}(y-x).
\]
Therefore, for a fixed $i \in \mathbb N$, we have that
\[
\lim_{x\rightarrow \infty} d^i(i+1,x) \leq \sum_{j=1}^\infty d^i(i+j,i+j+1) \leq \sum_{j=1}^\infty \frac{1}{1+j^2} < \infty.
\]
The finiteness of the negative part of $\mathbb R$ follows analogously.
Therefore $d^i$ is a finite metric (and therefore it isn't complete) for every $i \in \mathbb N$ despite $d^\infty$ is complete, what shows that the completeness of $d$ is not necessary to make $d^\infty$ intrinsic.

\end{example}

\begin{example}
\label{exemplo diverge}
Consider the additive group $\mathbb R$  endowed with the metric $d(x,y)=\vert \sqrt[3]{x} - \sqrt[3]{y} \vert$ and let $X=\{-1,0,\sqrt{2}\}$.
Suppose that the sequence of induced Hausdorff metric converges to a metric.
$\mathbb R$ doesn't admit a $C^0$-Carnot-Carath\'eodory-Finsler metric which isn't a $C^0$-Finsler metric, because the only distribution that can generate a $C^0$-Carnot-Carath\'eodory-Finsler metric is the whole tangent bundle.
Therefore the sequence converges to a $C^0$-Finsler metric.
But in this case it is straightforward to see that $\tilde F$ in (\ref{f tilde}) is infinite (just analyze small neighborhoods of the origin), what contradicts Theorem \ref{teorema principal}.
Therefore the sequence of induced Hausdorff metrics doesn't converge to a metric on $\mathbb R$.
\end{example}

\begin{example}
\label{flat torus}
Consider the flat torus $T=S^1 \times S^1$ with the canonical group operation.
We represent it as $[0,1]\times [0,1]$ endowed with the canonical metric with the opposite sides identified and endowed with the quotient metric $d$.
Consider $X=T-Q$, where $Q$ is the equivalence class of the square $(1/4,3/4) \times (1/4,3/4) \subset [0,1] \times [0,1]$ in $T^2$. 
It is not difficult to see that $d_X$ is the maximum metric when restricted to a small neighborhood of the identity element. 
Therefore $d^1$ is locally equal to the maximum metric. 
(Compare with Example 8.11 in \cite{BenettiFukuoka}).
In this case we have that $d^1 < d$, a relationship that doesn't happen in this work.
This example shows that the sequence of induced Hausdorff metrics are not only about increasing sequences of metrics.
\end{example}

\section{Final remarks}
\label{secao final}

\begin{remark}
We don't know if Theorem \ref{teorema principal} holds if $d$ isn't complete.
More in general, we don't know if the limit of an increasing sequence of intrinsic metrics is intrinsic.
\end{remark}

\begin{remark}
From a Lie group $G$ endowed with a metric satisfying mild conditions, we defined a discrete dynamical system $d,d^1, \ldots$ that converges to a right invariant $C^0$-Carnot-Carath\'eodory-Finsler metric.
This is the first example that shows that a sequence of metrics defined from the Hausdorff distance can have a regularizing effect.
Moreover, for a general compact subset $X \subset G$, the sequence of induced Hausdorff metrics is not necessarily increasing, as Example \ref{flat torus} shows, what makes the dynamics more interesting.
Due to these facts we think that it is worthwhile to study the properties of the sequence of induced Hausdorff metrics in other situations.
\end{remark}


\end{document}